\documentclass[12pt]{article}
\usepackage{amssymb}
\usepackage{amsmath}
\newtheorem{case}{\normalfont Case}
\newtheorem{example}{\normalfont\bf Example}
\newcommand{\sech}{\mathrm{sech}}
\newcommand{\csch}{\mathrm{csch}}
\newcommand{\sn}{\mathrm{sn}}
\newcommand{\cn}{\mathrm{cn}}
\newcommand{\dn}{\mathrm{dn}}
\newcommand{\ds}{\mathrm{ds}}
\newcommand{\ns}{\mathrm{ns}}
\newcommand{\cs}{\mathrm{cs}}

\begin{document}
\title{Variable separated ODE method--A powerful tool for testing
       traveling wave solutions of nonlinear equations}
\author{Sirendaoreji\\{\small Mathematical Science College,Inner Mongolia Normal
        University,}\\{\small Huhhot 010022,Inner Mongolia, PR China.}\\
        {\small Email:siren@imnu.edu.cn}}
\date{October,11,2018}
\maketitle
\noindent
{{\bf Abstract:} The variable separated ODE method is extended by choosing the
additional variable separated equation as the general elliptic equation.More exact
traveling wave solutions of nonlinear equations are obtained by using the method of
comparison of coefficients and the known solutions of the auxiliary equation.}

\section{Introduction}
\label{sec-1}
In 2002, we established a direct algebraic method for solving
sine--Gordon type and sinh--Gordon type equations\cite{sd1},
which was named  as the variable separated ODE method by Wazwaz\cite{sd2}.
Further applications and generalizations were considered by
Wazwaz\cite{sd2,sd3,sd4,sd5,sd6}, Fu\cite{sd7}, Xie\cite{sd8,sd9,sd10}
and others\cite{sd11,sd12,sd13}, and they successfully solved many
generalized sine--Gordon type,sinh-Gordon type,cosine--Gordon,
cosh--Gordon type,sine--cosine--Gordon type, sinh--cosh--Gordon
type and Liouville type equations,etc. These studies showed that
the variable separated ODE method works effectively if the
equation involves sine, cosine, hyperbolic sine, and hyperbolic
cosine functions.If it is not the case, can we use the variable
separated ODE method to solve the nonlinear equations? In other
words,can we solve the commonly used nonlinear equations by
means of the variable separated ODE method? The purpose of this
paper is to give a positive answer to this question. The key to
solve the problem is how to select the additional variable separated
equation
\begin{equation}
\label{eq:1}
u^\prime(\xi)={\frac {du}{d\xi}}=G(u),
\end{equation}
which satisfied by the wave transformation $u(x,t)=u(\xi),\xi=x-{\omega}t$.
In general,the selection of Eq.(\ref{eq:1}) is can be performed by many
choices,such as the Riccati equation,the Bernoulli equation and the
general elliptic equation,etc. In particular,more exact traveling wave
solutions for the given nonlinear equations are obtained by choosing
the general elliptic equation as our additional variable separated
equation.The advantage of our method is that it does not need to
solve the nonlinear equations and the exact traveling wave
solutions of the nonlinear equations can be obtained by
comparing the coefficients of the related polynomials and using the
known solutions of the additional variable separated Eq.(\ref{eq:1}).
Therefore,the method is simpler than the existing direct algebraic
methods.\par
The paper is organized as follows.In section \ref{sec-2},
we shall give the description of the method.The introduction of
the auxiliary equations and the list of their solutions are also
presented.Some examples are presented in the section (\ref{sec-3}).
Section (\ref{sec-4}) will give some discussions.
\section{Description of the method}
\label{sec-2}
As described in [\ref{eq:1}],we first make the wave transformation
$u(x,t)=u(\xi)$,$\xi=x-{\omega}t$ to carry out a given nonlinear PDE
into an equivalent ODE.Substituting Eq.(\ref{eq:1}) into the given
ODE yields a system of algebraic equations that can be solved to
determine the unknown parameters.In this paper,the additional
variable separated Eq.(\ref{eq:1}) will be taken as the general
elliptic equation \cite{sd14,sd15,sd16,sd17,sd18} of the form
\begin{equation}
\label{eq:2}
{F^\prime}^2(\xi)=c_0+c_1F(\xi)+c_2F^2(\xi)+c_3F^3(\xi)+c_4F^4(\xi),
\end{equation}
where $c_i~(i=0,1,2,3,4)$ are constants.Here we shall use a new
classification of solutions for Eq.(\ref{eq:2}) as the following
five cases.
\begin{case}
\label{case-1}
\normalfont
$c_0=c_1=0$.
\begin{align*}
&F_1(\xi)={\frac {2c_2}{\varepsilon\sqrt{\Delta}\cosh\left(\sqrt{c_2}\xi\right)-c_3}},\Delta>0,c_2>0,\\
&F_2(\xi)={\frac {2c_2}{\varepsilon\sqrt{-\Delta}\sinh\left(\sqrt{c_2}\xi\right)-c_3}},\Delta<0,c_2>0,\\
&F_{3a}(\xi)={\frac {2c_2}{\varepsilon\sqrt{\Delta}\cos\left(\sqrt{-c_2}\xi\right)-c_3}},\\
&\quad\,F_{3b}(\xi)={\frac {2c_2}{\varepsilon\sqrt{\Delta}\sin\left(\sqrt{-c_2}\xi\right)-c_3}},\Delta>0,c_2<0,\\
&F_4(\xi)=-{\frac {c_2}{c_3}}\left[1+\varepsilon\tanh\left({\frac {\sqrt{c_2}}2}\xi\right)\right],\Delta=0,c_2>0,\\
&F_5(\xi)=-{\frac {c_2}{c_3}}\left[1+\varepsilon\coth\left({\frac {\sqrt{c_2}}2}\xi\right)\right],\Delta=0,c_2>0,\\
&F_6(\xi)={\frac {\varepsilon}{\sqrt{c_4}\xi}},c_2=c_3=0,c_4>0,\\
&F_7(\xi)={\frac {4c_3}{c_3^2\xi^2-4c_4}},c_2=0,
\end{align*}
where $\Delta=c_3^2-4c_2c_4,\varepsilon=\pm{1}$.
\end{case}
\begin{case}
\label{case-2}
\normalfont
$c_3=c_4=0$.
\begin{align*}
&F_8(\xi)=-{\frac {c_1}{2c_2}}+{\frac {\varepsilon\sqrt{\delta}}{2c_2}}\cosh\left(\sqrt{c_2}\xi\right),\delta>0,c_2>0,\\
&F_9(\xi)=-{\frac {c_1}{2c_2}}+{\frac {\varepsilon\sqrt{-\delta}}{2c_2}}\sinh\left(\sqrt{c_2}\xi\right),\delta<0,c_2>0,\\
&F_{10a}(\xi)=-{\frac {c_1}{2c_2}}+{\frac {\varepsilon\sqrt{\delta}}{2c_2}}\cos\left(\sqrt{-c_2}\xi\right),\\
&\quad\,F_{10b}(\xi)=-{\frac {c_1}{2c_2}}+{\frac {\varepsilon\sqrt{\delta}}{2c_2}}\sin\left(\sqrt{-c_2}\xi\right),\delta>0,c_2<0,\\
&F_{11}(\xi)=-{\frac {c_1}{2c_2}}+e^{\varepsilon\sqrt{c_2}\xi},\delta=0,c_2>0,\\
&F_{12}(\xi)=\varepsilon\sqrt{c_0}\xi,c_1=c_2=0,\\
&F_{13}(\xi)=-{\frac {c_0}{c_1}}+{\frac {c_1}4}\xi^2,c_2=0,
\end{align*}
where $\delta=c_1^2-4c_0c_2,\varepsilon=\pm{1}$.
\end{case}
\begin{case}
\label{case-3}
\normalfont
$c_1=c_3=0$.
\begin{align*}
&F_{14}(\xi)=\varepsilon\sqrt{-\frac {c_2}{2c_4}}\tanh\left(\sqrt{-\frac {c_2}2}\xi\right),\Delta_1=0,c_2<0,c_4>0,\\
&F_{15}(\xi)=\varepsilon\sqrt{-\frac {c_2}{2c_4}}\coth\left(\sqrt{-\frac {c_2}2}\xi\right),\Delta_1=0,c_2<0,c_4>0,\\
&F_{16a}(\xi)=\varepsilon\sqrt{\frac {c_2}{2c_4}}\tan\left(\sqrt{\frac {c_2}2}\xi\right),\\
&\quad\,F_{16b}(\xi)=\varepsilon\sqrt{\frac {c_2}{2c_4}}\cot\left(\sqrt{\frac {c_2}2}\xi\right),\Delta_1=0,c_2>0,c_4>0;\\
&F_{17}(\xi)=\sqrt{\frac {-c_2m^2}{c_4(m^2+1)}}\sn\left(\sqrt{\frac {-c_2}{m^2+1}}\xi\right), c_0={\frac {c_2^2m^2}{c_4(m^2+1)^2}},c_2<0,c_4>0,\\
&F_{18}(\xi)=\sqrt{\frac {-c_2m^2}{c_4(2m^2-1)}}\cn\left(\sqrt{\frac {c_2}{2m^2-1}}\xi\right), c_0={\frac {c_2^2m^2(m^2-1)}{c_4(2m^2-1)^2}},\\
&\qquad\,c_2>0,c_4<0,\\
&F_{19}(\xi)=\sqrt{\frac {-c_2}{c_4(2-m^2)}}\dn\left(\sqrt{\frac {c_2}{2-m^2}}\xi\right), c_0={\frac {c_2^2(1-m^2)}{c_4(2-m^2)^2}},c_2>0,c_4<0,\\
&F_{20}(\xi)=\varepsilon\left(-{\frac {4c_0}{c_4}}\right)^{\frac 14}\ds\left((-4c_0c_4)^{\frac 14}\xi,{\frac {\sqrt{2}}2}\right),c_2=0,c_0c_4<0,\\
&F_{21}(\xi)=\varepsilon\left({\frac {c_0}{c_4}}\right)^{\frac 14}\left[\ns\left(2(c_0c_4)^{\frac 14}\xi,{\frac {\sqrt{2}}2}\right)
     +\cs\left(2(c_0c_4)^{\frac 14}\xi,{\frac {\sqrt{2}}2}\right)\right],\\
&\qquad\,c_2=0,c_0c_4>0,
\end{align*}
where $\Delta_1=c_2^2-4c_0c_4,\varepsilon=\pm{1}$.
\end{case}
\begin{case}
\label{case-4}
\normalfont
$c_2=c_4=0$.
\begin{align*}
F_{22}(\xi)=\wp\left({\frac {\sqrt{c_3}}2}\xi,g_2,g_3\right),g_2=-{\frac {4c_1}{c_3}},g_3=-{\frac {4c_0}{c_3}},c_3>0.
\end{align*}
\end{case}
\begin{case}
\label{case-5}
\normalfont
$c_0=0$.
\begin{align*}
&F_{23}(\xi)=-{\frac {8c_2\tanh^2\left(\sqrt{-\frac {c_2}{12}}\xi\right)}
   {3c_3\left(3+\tanh^2\left(\sqrt{-\frac {c_2}{12}}\xi\right)\right)}},\\
&F_{24}(\xi)=-{\frac {8c_2\coth^2\left(\sqrt{-\frac {c_2}{12}}\xi\right)}
   {3c_3\left(3+\coth^2\left(\sqrt{-\frac {c_2}{12}}\xi\right)\right)}},c_2<0,c_1={\frac {8c_2^2}{27c_3}},c_4={\frac {c_3^2}{4c_2}};\\
&F_{25}(\xi)={\frac {8c_2\tan^2\left(\sqrt{\frac {c_2}{12}}\xi\right)}
   {3c_3\left(3-\tan^2\left(\sqrt{\frac {c_2}{12}}\xi\right)\right)}},\\
&F_{26}(\xi)={\frac {8c_2\cot^2\left(\sqrt{\frac {c_2}{12}}\xi\right)}
   {3c_3\left(3-\cot^2\left(\sqrt{\frac {c_2}{12}}\xi\right)\right)}},c_2>0,c_1={\frac {8c_2^2}{27c_3}},c_4={\frac {c_3^2}{4c_2}};\\
&F_{27}(\xi)=-{\frac {c_3}{4c_4}}\left[1+\varepsilon{\sn}\left({\frac {c_3}
   {4m\sqrt{c_4}}}\xi\right)\right],\\
&F_{28}(\xi)=-{\frac {c_3}{4c_4}}\left[1+{\frac \varepsilon{m{\sn}\left({\frac {c_3}
   {4m\sqrt{c_4}}}\xi\right)}}\right],\\
   &\qquad\,c_4>0,c_1={\frac {c_3^3(m^2-1)}{32m^2c_4^2}},c_2={\frac {c_3^2(5m^2-1)}{16m^2c_4}};\\
&F_{29}(\xi)=-{\frac {c_3}{4c_4}}\left[1+\varepsilon{m\sn}\left({\frac {c_3}
   {4\sqrt{c_4}}}\xi\right)\right],\\
&F_{30}(\xi)=-{\frac {c_3}{4c_4}}\left[1+{\frac \varepsilon{{\sn}\left({\frac {c_3}
   {4\sqrt{c_4}}}\xi\right)}}\right],\\
   &\qquad\,c_4>0,c_1={\frac {c_3^3(1-m^2)}{32c_4^2}},c_2={\frac {c_3^2(5-m^2)}{16c_4}};\\
&F_{31}(\xi)=-{\frac {c_3}{4c_4}}\left[1+\varepsilon{\cn}\left({-\frac {c_3}
   {4m\sqrt{-c_4}}}\xi\right)\right],\\
&F_{32}(\xi)=-{\frac {c_3}{4c_4}}\left[1+{\frac {\varepsilon\sqrt{1-m^2}{\sn}\left({-\frac {c_3}
   {4m\sqrt{-c_4}}}\xi\right)}{{\dn}\left({-\frac {c_3}{4m\sqrt{-c_4}}}\xi\right)}}\right],\\
   &\qquad\,c_4<0,c_1={\frac {c_3^3}{32m^2c_4^2}},c_2={\frac {c_3^2(4m^2+1)}{16m^2c_4}};
\end{align*}
\begin{align*}
&F_{33}(\xi)=-{\frac {c_3}{4c_4}}\left[1+{\frac {\varepsilon}{\sqrt{1-m^2}}}{\dn}
\left({\frac {c_3}{4\sqrt{c_4(m^2-1)}}}\xi\right)\right],\\
&F_{34}(\xi)=-{\frac {c_3}{4c_4}}\left[1+{\frac {\varepsilon}{{\dn}\left({\frac {c_3}
   {4\sqrt{c_4(m^2-1)}}}\xi\right)}}\right],\\
   &\qquad\,c_4<0,c_1={\frac {c_3^3m^2}{32c_4^2(m^2-1)}},c_2={\frac {c_3^2(5m^2-4)}{16c_4(m^2-1)}};\\
&F_{35}(\xi)=-{\frac {c_3}{4c_4}}\left[1+{\frac {\varepsilon}{{\cn}\left({\frac {c_3}
   {4\sqrt{c_4(1-m^2)}}}\xi\right)}}\right],\\
&F_{36}(\xi)=-{\frac {c_3}{4c_4}}\left[1+{\frac {\varepsilon{\dn}\left({\frac {c_3}
   {4\sqrt{c_4(1-m^2)}}}\xi\right)}{\sqrt{1-m^2}{\cn}\left({\frac {c_3}{4\sqrt{c_4(1-m^2)}}}\xi\right)}}\right],\\
&\qquad\,c_4>0,c_1={\frac {c_3^3}{32c_4^2(1-m^2)}},c_2={\frac {c_3^2(4m^2-5)}{16c_4(m^2-1)}};\\
&F_{37}(\xi)=-{\frac {c_3}{4c_4}}\left[1+\varepsilon{\dn}\left({-\frac {c_3}
   {4\sqrt{-c_4}}}\xi\right)\right],\\
&F_{38}(\xi)=-{\frac {c_3}{4c_4}}\left[1+{\frac {\varepsilon\sqrt{1-m^2}}{{\dn}\left({-\frac {c_3}
   {4\sqrt{-c_4}}}\xi\right)}}\right],c_4<0,c_1={\frac {c_3^3m^2}{32c_4^2}},c_2={\frac {c_3^2(m^2+4)}{16c_4}},
\end{align*}
where $\varepsilon=\pm{1}$.
\end{case}
\section{Illustrative examples}
\label{sec-3}
In this section we shall give three illustrative examples
for using our method to solve nonlinear equations.
\begin{example}
\normalfont
The modified Benjamin--Bona--Mahony(mBBM) equation
\begin{align}
\label{eq:3}
u_t+u_x+u^2u_x+u_{xxt}=0.
\end{align}
\end{example}
Substituting the wave transformation $u(x,t)=u(\xi),\xi=x-{\omega}t$ into (\ref{eq:3}) and integrating once with respect to $\xi$ leads the following ODE
\begin{align*}
\left(1-{\omega}\right)u(\xi)+{\frac 13}u^3(\xi)-{\omega}u^{\prime\prime}(\xi)+B=0,
\end{align*}
which can be rewritten as
\begin{align}
\label{eq:4}
u^{\prime\prime}(\xi)={\frac {B}\omega}+{\frac {1-\omega}{\omega}}u(\xi)
  +{\frac 1{3\omega}}u^3(\xi),
\end{align}
where $B$ is the integration constant.\par
In accordance with the form of Eq. (\ref{eq:4}),we can take the Eq.(\ref{eq:1})
as the the general elliptic equation for unknown function $u(\xi)$ of the form
\begin{align*}
{u^\prime}^2(\xi)=G^2(u)=c_0+c_1u(\xi)+c_2u^2(\xi)+c_3u^3(\xi)+c_4u^4(\xi),
\end{align*}
which has the second form
\begin{align}
\label{eq:5}
u^{\prime\prime}(\xi)={\frac {c_1}2}+c_2u(\xi)+{\frac {3c_3}2}u^2(\xi)+2c_4u^3(\xi),
\end{align}
where $c_i~(i=0,1,2,3,4)$ are undetermined parameters.\par
Comparing the coefficients of $u^j~(j=0,1,2,3)$ in Eq.(\ref{eq:4})
and Eq.(\ref{eq:5}) yields
\begin{align}
\label{eq:6}
c_0=c_0,c_1={\frac {2B}{\omega}},c_2={\frac {1-\omega}{\omega}},c_3=0,c_4={\frac 1{6\omega}}.
\end{align}
When $c_0=B=0$,if we take (\ref{eq:6}) into the solutions of the
general elliptic equation (\ref{eq:2}) given by Case \ref{case-1},
then the exact traveling wave solutions of the mBBM equation are
obtained as follows
\begin{align*}
&u_1(x,t)=\varepsilon\sqrt{6(1-\omega)}\csch\left(\sqrt{\frac {1-\omega}{\omega}}(x-{\omega}t)\right),0<\omega<1,\\
&u_2(x,t)=\varepsilon\sqrt{6(\omega-1)}\sec\left(\sqrt{\frac {\omega-1}{\omega}}(x-{\omega}t)\right),\omega>1,\\
&u_3(x,t)=\varepsilon\sqrt{6(\omega-1)}\csc\left(\sqrt{\frac {\omega-1}{\omega}}(x-{\omega}t)\right),\omega>1,\\
&u_4(x,t)={\frac {\sqrt{6}\varepsilon}{x-t+\xi_0}},\omega=1.
\end{align*}
When $c_0\not=0,B=0$,taking (\ref{eq:6}) into the solutions of
the general elliptic equation (\ref{eq:2}) given by Case \ref{case-3},
we get the following exact traveling wave solutions of the mBBM equation
\begin{align*}
&u_5(x,t)=\varepsilon\sqrt{3(\omega-1)}\tanh\left(\sqrt{\frac {\omega-1}{2\omega}}(x-{\omega}t)\right),\omega>1,\\
&u_6(x,t)=\varepsilon\sqrt{3(\omega-1)}\coth\left(\sqrt{\frac {\omega-1}{2\omega}}(x-{\omega}t)\right),\omega>1,\\
&u_7(x,t)=\varepsilon\sqrt{3(1-\omega)}\tan\left(\sqrt{\frac {1-\omega}{2\omega}}(x-{\omega}t)\right),0<\omega<1,\\
&u_8(x,t)=\varepsilon\sqrt{3(1-\omega)}\cot\left(\sqrt{\frac {1-\omega}{2\omega}}(x-{\omega}t)\right),0<\omega<1,
\end{align*}
\begin{align*}
&u_9(x,t)=\sqrt{\frac {6(\omega-1)m^2}{m^2+1}}\sn\left(\sqrt{\frac {\omega-1}{\omega(m^2+1)}}(x-{\omega}t),m\right),\omega>1,\\
&u_{10}(x,t)=\varepsilon\left(-24{c_0}\right)^{\frac 14}\ds\left(\left(-{\frac 23}{c_0}\right)^{\frac 14}(x-t),{\frac {\sqrt{2}}{2}}\right),\omega=1,c_0<0,\\
&u_{11}(x,t)=\varepsilon\left(6{c_0}\right)^{\frac 14}\left[\ns\left(2\left({\frac {c_0}6}\right)^{\frac 14}(x-t),{\frac {\sqrt{2}}{2}}\right)\right.\\
&\qquad+\left.\cs\left(2\left({\frac {c_0}6}\right)^{\frac 14}(x-t),{\frac {\sqrt{2}}{2}}\right)\right],\omega=1,c_0>0.
\end{align*}
For solutions $u_j~(j=5,6,7,8)$,the condition
$\Delta_1=c_2^2-4c_0c_4=\left({\frac {1-\omega}{\omega}}\right)^2-{\frac {2c_0}{3\omega}}=0$
leads $c_0={\frac {3(1-\omega)^2}{2\omega}}$.For solution $u_9$,we have
$c_0={\frac {6(1-\omega)^2m^2(m^2-1)}{\omega(2m^2-1)^2}}$.
But the solutions $u_j~(j=5,6,7,8,9)$ do not require
any condition and they are true for any arbitrary constant $c_0\not=0$.
Therefore,it does not need to calculate the value of $c_0$.
\begin{example}
\normalfont
Consider the nonlinear Schr\"odinger(NLS) equation
\begin{align}
\label{eq:7}
iu_t+{\alpha}u_{xx}+\beta|u|^2u=0,
\end{align}
where $\alpha$ and $\beta$ are constants.
\end{example}
Taking the transformation
\begin{align}
\label{eq:8}
u(x,t)=v(\xi)e^{i\eta},\xi=x-{\omega}t,\eta=kx+ct,
\end{align}
with undetermined constants $c,k$ and $\omega$ into Eq. (\ref{eq:7})
and separating the real part and the imaginary part yields
\begin{align*}
\left\{\begin{aligned}
&{-\omega}v^\prime(\xi)+2\alpha{k}v^\prime(\xi)=0,\\
&-cv(\xi)+\alpha\left(v^{\prime\prime}(\xi)-k^2v(\xi)\right)+{\beta}v^3(\xi)=0.
\end{aligned}\right.
\end{align*}
The first equation leads $k={\frac {\omega}{2\alpha}}$ and the second
equation becomes
\begin{align*}
-cv(\xi)+\alpha\left(v^{\prime\prime}(\xi)-{\frac {\omega^2}{4\alpha^2}}v(\xi)\right)+{\beta}v^3(\xi)=0.
\end{align*}
This equation can be rewritten as
\begin{align}
\label{eq:9}
v^{\prime\prime}(\xi)={\frac {\omega^2+4{\alpha}c}{4\alpha^2}}v(\xi)-{\frac {\beta}{\alpha}}v^3(\xi).
\end{align}
If we choose the Eq. (\ref{eq:1}) as the general elliptic equation
\begin{align*}
{v^\prime}^2(\xi)=G^2(v)=c_0+c_1v(\xi)+c_2v^2(\xi)+c_3v^3(\xi)+c_4v^4(\xi),
\end{align*}
then by differentiating the above equation with respect to $\xi$,we obtain
\begin{align}
\label{eq:10}
v^{\prime\prime}(\xi)={\frac {c_1}2}+c_2v(\xi)+{\frac {3c_3}2}v^2(\xi)+2c_4v^3(\xi).
\end{align}
By comparing the coefficients of $v^j~(j=0,1,2,3)$ in Eq. (\ref{eq:9})
and Eq. (\ref{eq:10}),the unknown parameters are now determined to be
\begin{align}
\label{eq:11}
c_0=c_0,c_1=0,c_2={\frac {\omega^2+4\alpha{c}}{4\alpha^2}},c_3=0,c_4=-{\frac {\beta}{2\alpha}}.
\end{align}
When $c_0=0$,inserting (\ref{eq:11}) into the solutions of
the general elliptic equation given by Case \ref{case-1}
and using (\ref{eq:8}) we get the following exact traveling
wave solutions of NLS equation
\begin{align*}
&u_1(x,t)=\varepsilon\sqrt{\frac {\omega^2+4{\alpha}c}{2\alpha\beta}}
  \sech\left[{\frac 12}\sqrt{\frac {\omega^2+4{\alpha}c}{\alpha^2}}\left(x-{\omega}t\right)\right]
  e^{i\left({\frac {\omega}{2\alpha}}x+ct\right)},\\
  &\qquad\omega^2+4{\alpha}c>0,\alpha\beta>0,\\
&u_2(x,t)=\varepsilon\sqrt{-\frac {\omega^2+4{\alpha}c}{2\alpha\beta}}
  \csch\left[{\frac 12}\sqrt{\frac {\omega^2+4{\alpha}c}{\alpha^2}}\left(x-{\omega}t\right)\right]
  e^{i\left({\frac {\omega}{2\alpha}}x+ct\right)},\\
  &\qquad\omega^2+4{\alpha}c>0,\alpha\beta<0,\\
&u_3(x,t)=\varepsilon\sqrt{\frac {\omega^2+4{\alpha}c}{2\alpha\beta}}
  \sec\left[{\frac 12}\sqrt{-\frac {\omega^2+4{\alpha}c}{\alpha^2}}\left(x-{\omega}t\right)\right]
  e^{i\left({\frac {\omega}{2\alpha}}x+ct\right)},\\
  &\qquad\omega^2+4{\alpha}c<0,\alpha\beta<0,\\
&u_4(x,t)=\varepsilon\sqrt{\frac {\omega^2+4{\alpha}c}{2\alpha\beta}}
  \csc\left[{\frac 12}\sqrt{-\frac {\omega^2+4{\alpha}c}{\alpha^2}}\left(x-{\omega}t\right)\right]
  e^{i\left({\frac {\omega}{2\alpha}}x+ct\right)},\\
  &\qquad\omega^2+4{\alpha}c<0,\alpha\beta<0,\\
&u_5(x,t)={\frac {\varepsilon}{\sqrt{-\frac {\beta}{2\alpha}}(x-{\omega}t+\xi_0)}}
  e^{i\left({\frac {\omega}{2\alpha}}x+ct\right)},\omega^2+4{\alpha}c=0,\alpha\beta>0,
\end{align*}
where $\omega,\xi_0$ are arbitrary constants.\par
When $c_0\not=0$,substituting (\ref{eq:11}) into the solutions
of general elliptic equation given by Case \ref{case-3},we obtain
the exact traveling wave solutions of NLS equation as follows
\begin{align*}
&u_6(x,t)={\frac {\varepsilon}2}\sqrt{\frac {\omega^2+4{\alpha}c}{\alpha\beta}}
  \tanh\left[{\frac 14}\sqrt{-\frac {2(\omega^2+4{\alpha}c)}{\alpha^2}}\left(x-{\omega}t\right)\right]
  e^{i\left({\frac {\omega}{2\alpha}}x+ct\right)},\\
  &\qquad\omega^2+4{\alpha}c<0,\alpha\beta<0,\\
&u_7(x,t)={\frac {\varepsilon}2}\sqrt{\frac {\omega^2+4{\alpha}c}{\alpha\beta}}
  \coth\left[{\frac 14}\sqrt{-\frac {2(\omega^2+4{\alpha}c)}{\alpha^2}}\left(x-{\omega}t\right)\right]
  e^{i\left({\frac {\omega}{2\alpha}}x+ct\right)},\\
  &\qquad\omega^2+4{\alpha}c<0,\alpha\beta<0,
\end{align*}
\begin{align*}
&u_8(x,t)={\frac {\varepsilon}2}\sqrt{-\frac {\omega^2+4{\alpha}c}{\alpha\beta}}
  \tan\left[{\frac 14}\sqrt{\frac {2(\omega^2+4{\alpha}c)}{\alpha^2}}\left(x-{\omega}t\right)\right]
  e^{i\left({\frac {\omega}{2\alpha}}x+ct\right)},\\
  &\qquad\omega^2+4{\alpha}c>0,\alpha\beta<0,\\
&u_9(x,t)={\frac {\varepsilon}2}\sqrt{-\frac {\omega^2+4{\alpha}c}{\alpha\beta}}
  \cot\left[{\frac 14}\sqrt{\frac {2(\omega^2+4{\alpha}c)}{\alpha^2}}\left(x-{\omega}t\right)\right]
  e^{i\left({\frac {\omega}{2\alpha}}x+ct\right)},\\
  &\qquad\omega^2+4{\alpha}c>0,\alpha\beta<0,\\
&u_{10}(x,t)=\sqrt{\frac {m^2(\omega^2+4{\alpha}c)}{2\alpha\beta(2m^2-1)}}
  \cn\left({\frac 12}\sqrt{\frac {\omega^2+4{\alpha}c}{\alpha^2(2m^2-1)}}\left(x-{\omega}t\right),m\right)
  e^{i\left({\frac {\omega}{2\alpha}}x+ct\right)},\\
  &\qquad\omega^2+4{\alpha}c>0,\alpha\beta>0,{\frac 12}<m^2<1,\\
&u_{11}(x,t)=\sqrt{\frac {m^2(\omega^2+4{\alpha}c)}{2\alpha\beta(m^2+1)}}
  \sn\left({\frac 12}\sqrt{-\frac {\omega^2+4{\alpha}c}{\alpha^2(m^2+1)}}\left(x-{\omega}t\right),m\right)
  e^{i\left({\frac {\omega}{2\alpha}}x+ct\right)},\\
  &\qquad\omega^2+4{\alpha}c<0,\alpha\beta<0,\\
&u_{12}(x,t)=\sqrt{\frac {\omega^2+4{\alpha}c}{2\alpha\beta(2-m^2)}}
  \dn\left({\frac 12}\sqrt{\frac {\omega^2+4{\alpha}c}{\alpha^2(2-m^2)}}\left(x-{\omega}t\right),m\right)
  e^{i\left({\frac {\omega}{2\alpha}}x+ct\right)},\\
&\qquad\omega^2+4{\alpha}c>0,\alpha\beta>0,\\
&u_{13}(x,t)=\varepsilon\left({-\frac {2\omega^2c_0}{\beta{c}}}\right)^{\frac 14}
  \ds\left(\left(-{\frac {8\beta{c_0}c}{\omega^2}}\right)^{\frac 14}\left(x-{\omega}t\right),{\frac {\sqrt{2}}{2}}\right)
  e^{-i\left({\frac {2c}{\omega}}x-ct\right)},\\
  &\qquad{\omega^2+4\alpha{c}}=0,\beta{c_0}c<0,\\
&u_{14}(x,t)=\varepsilon\left({\frac {\omega^2c_0}{2\beta{c}}}\right)^{\frac 14}\left[
  \ns\left(2\left({\frac {2\beta{c_0}c}{\omega^2}}\right)^{\frac 14}\left(x-{\omega}t\right),{\frac {\sqrt{2}}{2}}\right)\right.\\
  &\qquad\left.+\cs\left(2\left({\frac {2\beta{c_0}c}{\omega^2}}\right)^{\frac 14}\left(x-{\omega}t\right),
  {\frac {\sqrt{2}}{2}}\right)\right]e^{-i\left({\frac {2c}{\omega}}x-ct\right)},
  \omega^2+4\alpha{c}=0,\beta{c_0}c>0.
\end{align*}
It is noted that the solutions $u_j~(j=10,11,12)$ do not
require any condition,so it does not need to calculate $c_0$
from the conditions of solutions.
\begin{example}
\normalfont
The combined KdV--mKdV equation
\begin{eqnarray}
\label{eq:12}
u_t+6\left(\alpha{u}+\beta{u^2}\right)u_x+\gamma{u_{xxx}}=0,
\end{eqnarray}
in which $\alpha,\beta,\gamma$ are constants.
\end{example}\par
Substituting the wave transformation $u(x,t)=u(\xi),\xi=x-{\omega}t$
into Eq. (\ref{eq:12}) and integrating once with respect to $\xi$,
we obtain
\begin{eqnarray*}
-\omega{u}+3\alpha{u^2}+2\beta{u^3}+{\gamma}u^{\prime\prime}-C=0,
\end{eqnarray*}
where $C$ is the constant of integration.Now this equation can
be rewritten as
\begin{equation}
\label{eq:13}
u^{\prime\prime}={\frac C\gamma}+{\frac {\omega}{\gamma}}u-{\frac {3\alpha}{\gamma}}u^2
-{\frac {2\beta}{\gamma}}u^3.
\end{equation}
Suppose that the $u(\xi)$ satisfies the general elliptic equation,
then we have
\begin{equation}
\label{eq:14}
u^{\prime\prime}={\frac {c_1}2}+c_2u+{\frac {3c_3}2}u^2+2c_4u^3.
\end{equation}
Equating the coefficients of $u^j~(j=0,1,2,3)$ in Eq. (\ref{eq:13})
and Eq. (\ref{eq:14}) yields
\begin{equation}
\label{eq:15}
c_0=c_0,c_1={\frac {2C}{\gamma}},c_2={\frac {\omega}{\gamma}},
c_3={\frac {2\alpha}{\gamma}},c_4=-{\frac {\beta}{\gamma}}.
\end{equation}
\par\noindent
When $c_0=c_1=0$,inserting (\ref{eq:15}) into the solutions
in Case \ref{case-1} will give the following exact traveling
wave solutions of the combined KdV--mKdV equation
\begin{align*}
&u_1(x,t)={\frac {\omega}{\varepsilon\sqrt{\alpha^2+\beta\omega}\cosh\sqrt{\frac {\omega}{\gamma}}\left(x-{\omega}t\right)+\alpha}},\alpha^2+\beta\omega>0,\omega\gamma>0,\\
&u_2(x,t)={\frac {\omega}{\varepsilon\sqrt{-\left(\alpha^2+\beta\omega\right)}\sinh\sqrt{\frac {\omega}{\gamma}}\left(x-{\omega}t\right)+\alpha}},\alpha^2+\beta\omega<0,\omega\gamma>0,\\
&u_3(x,t)={\frac {\omega}{\varepsilon\sqrt{\alpha^2+\beta\omega}\cos\sqrt{-\frac {\omega}{\gamma}}\left(x-{\omega}t\right)+\alpha}},\alpha^2+\beta\omega>0,\omega\gamma<0,\\
&u_4(x,t)={\frac {\omega}{\varepsilon\sqrt{\alpha^2+\beta\omega}\sin\sqrt{-\frac {\omega}{\gamma}}\left(x-{\omega}t\right)+\alpha}},\alpha^2+\beta\omega>0,\omega\gamma<0,\\
&u_5(x,t)=-{\frac {\alpha}{2\beta}}\left[1+\varepsilon\tanh{\frac 12}\sqrt{-\frac {\alpha^2}{\beta\gamma}}\left(x+{\frac {\alpha^2}{\beta}}t\right)\right],\beta\gamma<0,\\
&u_6(x,t)=-{\frac {\alpha}{2\beta}}\left[1+\varepsilon\coth{\frac 12}\sqrt{-\frac {\alpha^2}{\beta\gamma}}\left(x+{\frac {\alpha^2}{\beta}}t\right)\right],\beta\gamma<0,
\end{align*}
and the stationary solution
\begin{align*}
&u_7(x,t)=-{\frac {2\alpha\gamma}{\alpha^2x^2+\beta\gamma}},\omega=0.
\end{align*}
When $c_0=0$,we need to consider the following seven cases.\par\noindent
(1).\,If $c_1={\frac {8c_2^2}{27c_3}},c_4={\frac {c_3^2}{4c_2}}$
then we obtain from (\ref{eq:15}) that
\begin{equation}
\label{eq:16}
c_1=-{\frac {4\alpha^3}{27\beta^2\gamma}},c_2=-{\frac {\alpha^2}{\beta\gamma}},
c_3=-{\frac {2\alpha}{\gamma}},\omega=-{\frac {\alpha^2}{\gamma}},C=-{\frac {2\alpha^3}{27\beta^2}}.
\end{equation}
To take (\ref{eq:16}) into $F_j~(j=23,24,25,26)$,we obtain the following
exact traveling solutions of the combined KdV--mKdV equation
\begin{align*}
&u_8(x,t)=-{\frac {4\alpha\tanh^2{\frac 16}\sqrt{\frac {3\alpha^2}{\beta\gamma}}\left(x+{\frac {\alpha^2}{\beta}}t\right)}{3\beta\left(3+\tanh^2{\frac 16}\sqrt{\frac {3\alpha^2}{\beta\gamma}}\left(x+{\frac {\alpha^2}{\beta}}t\right)\right)^2}},\beta\gamma>0,
\end{align*}
\begin{align*}
&u_9(x,t)=-{\frac {4\alpha\coth^2{\frac 16}\sqrt{\frac {3\alpha^2}{\beta\gamma}}\left(x+{\frac {\alpha^2}{\beta}}t\right)}{3\beta\left(3+\coth^2{\frac 16}\sqrt{\frac {3\alpha^2}{\beta\gamma}}\left(x+{\frac {\alpha^2}{\beta}}t\right)\right)^2}},\beta\gamma>0,\\
&u_{10}(x,t)={\frac {4\alpha\tan^2{\frac 16}\sqrt{-\frac {3\alpha^2}{\beta\gamma}}\left(x+{\frac {\alpha^2}{\beta}}t\right)}{3\beta\left(3-\tan^2{\frac 16}\sqrt{-\frac {3\alpha^2}{\beta\gamma}}\left(x+{\frac {\alpha^2}{\beta}}t\right)\right)^2}},\beta\gamma<0,\\
&u_{11}(x,t)={\frac {4\alpha\cot^2{\frac 16}\sqrt{-\frac {3\alpha^2}{\beta\gamma}}\left(x+{\frac {\alpha^2}{\beta}}t\right)}{3\beta\left(3-\cot^2{\frac 16}\sqrt{-\frac {3\alpha^2}{\beta\gamma}}\left(x+{\frac {\alpha^2}{\beta}}t\right)\right)^2}},\beta\gamma<0.
\end{align*}
\par\noindent
(2).\,If $c_1={\frac {c_3^3(m^2-1)}{32m^2c_4^2}},c_2={\frac {c_3^2(5m^2-1)}{16m^2c_4}}$,
then (\ref{eq:15}) gives the following parameters
\begin{align}
\label{eq:17}
&c_1=-{\frac {\alpha^3(m^2-1)}{4m^2\beta^2\gamma}},
c_2=-{\frac {\alpha^2(5m^2-1)}{4m^2\beta\gamma}},
c_3=-{\frac {2\alpha}{\gamma}},\nonumber\\
&\quad\,c_4=-{\frac {\beta}{\gamma}},
\omega=-{\frac {\alpha^2(5m^2-1)}{4m^2\beta}},
C=-{\frac {\alpha^3(m^2-1)}{8m^2\beta^2}}.
\end{align}
Inserting (\ref{eq:17}) into $F_j~(j=27,28)$ yields the exact
solutions of the combined KdV--mKdV equation as follows
\begin{align*}
&u_{12}(x,t)=-{\frac {\alpha}{2\beta}}\left[1+\varepsilon\sn{\frac {\alpha}{2m\beta}}\sqrt{-\frac {\beta}{\gamma}}\left(x+{\frac {\alpha^2(5m^2-1)}{4m^2\beta}}t\right)\right],\beta\gamma<0,\\
&u_{13}(x,t)=-{\frac {\alpha}{2\beta}}\left[1+{\frac {\varepsilon}{m\sn{\frac {\alpha}{2m\beta}}\sqrt{-\frac {\beta}{\gamma}}\left(x+{\frac {\alpha^2(5m^2-1)}{4m^2\beta}}t\right)}}\right],\beta\gamma<0.
\end{align*}
\par\noindent
(3).\,If $c_1={\frac {c_3^3(1-m^2)}{32c_4^2}},c_2={\frac {c_3^2(5-m^2)}{16c_4}}$,then
it solves from (\ref{eq:15}) that
\begin{align}
\label{eq:17}
&c_1=-{\frac {\alpha^3(1-m^2)}{4\beta^2\gamma}},
c_2=-{\frac {\alpha^2(5-m^2)}{4\beta\gamma}},
c_3=-{\frac {2\alpha}{\gamma}},\nonumber\\
&c_4=-{\frac {\beta}{\gamma}},\omega={\frac {\alpha^2(m^2-5)}{\beta}},
C=-{\frac {\alpha^3(m^2-1)}{8m^2\beta^2}}.
\end{align}
Substituting (\ref{eq:17}) into $F_j~(j=29,30)$,we obtain the exact
traveling wave solutions of Eq. (\ref{eq:12}) as following
\begin{align*}
&u_{14}(x,t)=-{\frac {\alpha}{2\beta}}\left[1-\varepsilon{m}\sn{\frac {\alpha}{2\gamma}}\sqrt{-\frac {\gamma}{\beta}}\left(x-{\frac {\alpha^2(m^2-5)}{\beta}}t\right)\right],\beta\gamma<0,\\
&u_{15}(x,t)=-{\frac {\alpha}{2\beta}}\left[1-{\frac {\varepsilon}{\sn{\frac {\alpha}{2\gamma}}\sqrt{-\frac {\gamma}{\beta}}\left(x-{\frac {\alpha^2(m^2-5)}{\beta}}t\right)}}\right],\beta\gamma<0.
\end{align*}
\par\noindent
(4).\,If $c_1={\frac {c_3^3}{32m^2c_4^2}},c_2={\frac {c_3^2(4m^2+1)}{16m^2c_4}}$,
then we obtain from (\ref{eq:15}) that
\begin{align}
\label{eq:18}
&c_1=-{\frac {\alpha^3}{4m^2\beta^2\gamma}},
c_2=-{\frac {\alpha^2(4m^2+1)}{4m^2\beta\gamma}},
c_3=-{\frac {2\alpha}{\gamma}},\nonumber\\
&c_4=-{\frac {\beta}{\gamma}},\omega=-{\frac {\alpha^2(4m^2+1)}{4m^2\beta}},
C=-{\frac {\alpha^3}{8m^2\beta^2}}.
\end{align}
To take (\ref{eq:18}) into $F_j~(j=31,32)$ we obtain the following
two types of exact traveling wave solutions for Eq. (\ref{eq:12})
\begin{align*}
&u_{16}(x,t)=-{\frac {\alpha}{2\beta}}\left[1+\varepsilon\cn{\frac {\alpha}{2m\gamma}}\sqrt{\frac {\gamma}{\beta}}\left(x+{\frac {\alpha^2(4m^2+1)}{4m^2\beta}}t\right)\right],\beta\gamma>0,\\
&u_{17}(x,t)=-{\frac {\alpha}{2\beta}}\left[1+{\frac {\varepsilon\sqrt{1-m^2}\sn{\frac {\alpha}{2m\gamma}}\sqrt{\frac {\gamma}{\beta}}\left(x+{\frac {\alpha^2(4m^2+1)}{4m^2\beta}}t\right)}{\dn{\frac {\alpha}{2m\gamma}}\sqrt{\frac {\gamma}{\beta}}\left(x+{\frac {\alpha^2(4m^2+1)}{4m^2\beta}}t\right)}}\right],\beta\gamma>0.
\end{align*}
\par\noindent
(5).\,If $c_1={\frac {c_3^3m^2}{32c_4^2(m^2-1)}},c_2={\frac {c_3^2(5m^2-4)}{16c_4(m^2-1)}}$,
then (\ref{eq:15}) is solved as
\begin{align}
\label{eq:19}
&c_1=-{\frac {\alpha^3m^2}{4\beta^2\gamma(m^2-1)}},
c_2=-{\frac {\alpha^2(4m^2+1)}{4m^2\beta\gamma}},
c_3=-{\frac {2\alpha}{\gamma}},\nonumber\\
&c_4=-{\frac {\beta}{\gamma}},\omega=-{\frac {\alpha^2(5m^2-4)}{4\beta(m^2-1)}},
C=-{\frac {\alpha^3m^2}{8\beta^2(m^2-1)}}.
\end{align}
Inserting (\ref{eq:19}) into $F_j~(j=33,34)$ leads the exact traveling
wave solutions of Eq. (\ref{eq:12})
\begin{align*}
&u_{18}(x,t)=-{\frac {\alpha}{2\beta}}\left[1+{\frac {\varepsilon}{\sqrt{1-m^2}}}\dn{\frac {\alpha}{2\gamma}}\sqrt{\frac {\gamma}{\beta(1-m^2)}}\left(x+{\frac {\alpha^2(5m^2-4)}{4\beta(m^2-1)}}t\right)\right],\\
&u_{19}(x,t)=-{\frac {\alpha}{2\beta}}\left[1+{\frac {\varepsilon}{\dn{\frac {\alpha}{2\gamma}}\sqrt{\frac {\gamma}{\beta(1-m^2)}}\left(x+{\frac {\alpha^2(5m^2-4)}{4\beta(m^2-1)}}t\right)}}\right],
\end{align*}
where $\beta\gamma>0$.
\par\noindent
(6).\,If $c_1={\frac {c_3^3}{32c_4^2(1-m^2)}},c_2={\frac {c_3^2(4m^2-5)}{16c_4(m^2-1)}}$,then
(\ref{eq:15}) gives the following parameters
\begin{align}
\label{eq:20}
&c_1={\frac {\alpha^3}{4\beta^2\gamma(m^2-1)}},
c_2=-{\frac {\alpha^2(4m^2+1)}{4m^2\beta\gamma}},
c_3=-{\frac {2\alpha}{\gamma}},\nonumber\\
&c_4=-{\frac {\beta}{\gamma}},
\omega=-{\frac {\alpha^2(4m^2-5)}{4\beta(m^2-1)}},
C={\frac {\alpha^3}{8\beta^2(m^2-1)}}.
\end{align}
Substituting (\ref{eq:20}) into $F_j~(j=35,36)$ we obtain the
exact traveling wave solutions of Eq. (\ref{eq:12}) as following
\begin{align*}
&u_{20}(x,t)=-{\frac {\alpha}{2\beta}}\left[1+{\frac {\varepsilon}{\cn{\frac {\alpha}{2\gamma}}\sqrt{-\frac {\gamma}{\beta(1-m^2)}}\left(x+{\frac {\alpha^2(4m^2-5)}{4\beta(m^2-1)}}t\right)}}\right],\beta\gamma<0,\\
&u_{21}(x,t)=-{\frac {\alpha}{2\beta}}\left[1-{\frac {\varepsilon\dn{\frac {\alpha}{2\gamma}}\sqrt{-\frac {\gamma}{\beta(1-m^2)}}\left(x+{\frac {\alpha^2(4m^2-5)}{4\beta(m^2-1)}}t\right)}{\sqrt{1-m^2}\sn{\frac {\alpha}{2\gamma}}\sqrt{-\frac {\gamma}{\beta(1-m^2)}}\left(x+{\frac {\alpha^2(4m^2-5)}{4\beta(m^2-1)}}t\right)}}\right],\beta\gamma<0,
\end{align*}
\par\noindent
(7).\,If $c_1={\frac {c_3^3m^2}{32c_4^2}},c_2={\frac {c_3^2(m^2+4)}{16c_4}}$,
then (\ref{eq:15}) leads the following parameters
\begin{align}
\label{eq:21}
&c_1=-{\frac {\alpha^3m^2}{4\beta^2\gamma}},
c_2=-{\frac {\alpha^2(4m^2+1)}{4m^2\beta\gamma}},
c_3=-{\frac {2\alpha}{\gamma}},\nonumber\\
&c_4=-{\frac {\beta}{\gamma}},
\omega=-{\frac {\alpha^2(m^2+4)}{4\beta}},
C=-{\frac {\alpha^3m^2}{8\beta^2}}.
\end{align}
Taking (\ref{eq:21}) into $F_j~(j=37,38)$ we get the following exact traveling
wave solutions to Eq. (\ref{eq:12})
\begin{align*}
&u_{22}(x,t)=-{\frac {\alpha}{2\beta}}\left[1+\varepsilon\dn{\frac {\alpha}{2\gamma}}
\sqrt{\frac {\gamma}{\beta}}\left(x+{\frac {\alpha^2(m^2+4)}{4\beta}}t\right)\right],\beta\gamma>0,\\
&u_{23}(x,t)=-{\frac {\alpha}{2\beta}}\left[1+{\frac {\varepsilon\sqrt{1-m^2}}{\dn{\frac {\alpha}{2\gamma}}
\sqrt{\frac {\gamma}{\beta}}\left(x+{\frac {\alpha^2(m^2+4)}{4\beta}}t\right)}}\right],\beta\gamma>0.
\end{align*}
\section{Discussions}
\label{sec-4}
By our method one can use the known solutions of the additional
variable separated equation Eq. (\ref{eq:1}) to test the traveling
wave solutions of a given nonlinear equation.This is convenient
for judging the solutions of a given nonlinear equation without
solving the nonlinear equation.The method also does not require
the expansion expression of the traveling wave solution, so it
does not need to determine the order of the traveling wave
solutions from the homogeneous balancing principle.In addition,
we can neglect some conditions of Eq. (\ref{eq:1}) and determine
the integration constants by using the elementary comparison
method.Therefore, the method can be used to find  more
exact solutions of wide class of nonlinear equations.
\section*{Acknowledgement}
This work has been supported by the National Natural Science Foundation of
China under Grant No.11361040.

\bibliography{sdbibfile}

\end{document}